\documentclass[12pt]{article}

\usepackage{amsfonts,amssymb,amsthm,amsmath,latexsym}
\usepackage{subeqnarray,eqsection,indent,url}
\usepackage{bm}
\usepackage{graphicx}

\date{June 6, 2011 \\[1mm] revised January 21, 2012}   %% REMEMBER TO PUT FINAL DATE HERE!!!!!

\begin{document}

\title{\vspace*{-1cm}
   The leading root of the partial theta function
   % \\
   % and some generalizations(???)
   %% \\[4mm]
   %% $\Theta_0(x,y) = \displaystyle \sum\limits_{n=0}^\infty x^n y^{n(n-1)/2}$
      }

\author{
     \\
     {\small Alan D.~Sokal\thanks{Also at Department of Mathematics,
           University College London, London WC1E 6BT, England.}}  \\[-2mm]
     {\small\it Department of Physics}       \\[-2mm]
     {\small\it New York University}         \\[-2mm]
     {\small\it 4 Washington Place}          \\[-2mm]
     {\small\it New York, NY 10003 USA}      \\[-2mm]
     {\small\tt sokal@nyu.edu}               \\[-2mm]
     {\protect\makebox[5in]{\quad}}  % To force authors' names to be written
                                     %   vertically, one above another.
                                     % (\author seems to put them side-by-side
                                     %   if there is room.)
     \\
}

\maketitle
\thispagestyle{empty}   % Suppress page number on front page.

\begin{abstract}
I study the leading root $x_0(y)$ of the partial theta function
$\Theta_0(x,y) = \sum_{n=0}^\infty x^n y^{n(n-1)/2}$,
considered as a formal power series.
I prove that all the coefficients of $-x_0(y)$ are strictly positive.
Indeed, I prove the stronger results that all the coefficients of $-1/x_0(y)$
after the constant term 1 are strictly negative,
and all the coefficients of $1/x_0(y)^2$
after the constant term 1 are strictly negative
except for the vanishing coefficient of $y^3$.
\end{abstract}

\bigskip
\noindent
{\bf Key Words:}  Partial theta function, Rogers--Ramanujan function,
$q$-series, formal power series, root, implicit function theorem.

\bigskip
\noindent
{\bf Mathematics Subject Classification (MSC 2000) codes:}
05A15 (Primary);
05A19, 05A20, 05A30, 05C30, 11B65, 11P84, 30D20, 33D15, 33D65 (Secondary).

\clearpage

\newtheorem{theorem}{Theorem}[section]
\newtheorem{proposition}[theorem]{Proposition}
\newtheorem{lemma}[theorem]{Lemma}
\newtheorem{corollary}[theorem]{Corollary}
\newtheorem{definition}[theorem]{Definition}
\newtheorem{conjecture}[theorem]{Conjecture}
\newtheorem{question}[theorem]{Question}
\newtheorem{example}[theorem]{Example}

\renewcommand{\theenumi}{\alph{enumi}}
\renewcommand{\labelenumi}{(\theenumi)}
\def\eop{\hbox{\kern1pt\vrule height6pt width4pt
depth1pt\kern1pt}\medskip}
\def\prf{\par\noindent{\bf Proof.\enspace}\rm}
\def\rmk{\par\medskip\noindent{\bf Remark\enspace}\rm}

\newcommand{\be}{\begin{equation}}
\newcommand{\ee}{\end{equation}}
\newcommand{\<}{\langle}
\renewcommand{\>}{\rangle}
\newcommand{\widebar}{\overline}
\def\reff#1{(\protect\ref{#1})}
\def\spose#1{\hbox to 0pt{#1\hss}}
\def\ltapprox{\mathrel{\spose{\lower 3pt\hbox{$\mathchar"218$}}
    \raise 2.0pt\hbox{$\mathchar"13C$}}}
\def\gtapprox{\mathrel{\spose{\lower 3pt\hbox{$\mathchar"218$}}
    \raise 2.0pt\hbox{$\mathchar"13E$}}}
\def\textprime{${}^\prime$}
\def\proof{\par\medskip\noindent{\sc Proof.\ }}
\def\firstproof{\par\medskip\noindent{\sc First Proof.\ }}
\def\secondproof{\par\medskip\noindent{\sc Second Proof.\ }}
\def\qed{ $\square$ \bigskip}
\def\proofof#1{\bigskip\noindent{\sc Proof of #1.\ }}
\def\half{ {1 \over 2} }
\def\third{ {1 \over 3} }
\def\twothird{ {2 \over 3} }
\def\smfrac#1#2{{\textstyle{#1\over #2}}}
\def\smhalf{ {\smfrac{1}{2}} }
\newcommand{\real}{\mathop{\rm Re}\nolimits}
\renewcommand{\Re}{\mathop{\rm Re}\nolimits}
\newcommand{\imag}{\mathop{\rm Im}\nolimits}
\renewcommand{\Im}{\mathop{\rm Im}\nolimits}
\newcommand{\sgn}{\mathop{\rm sgn}\nolimits}
\newcommand{\tr}{\mathop{\rm tr}\nolimits}
\newcommand{\supp}{\mathop{\rm supp}\nolimits}
\def\hboxscript#1{ {\hbox{\scriptsize\em #1}} }
\renewcommand{\emptyset}{\varnothing}

\newcommand{\restrict}{\upharpoonright}

\newcommand{\scra}{{\mathcal{A}}}
\newcommand{\scrb}{{\mathcal{B}}}
\newcommand{\scrc}{{\mathcal{C}}}
\newcommand{\scrd}{{\mathcal{D}}}
\newcommand{\scre}{{\mathcal{E}}}
\newcommand{\scrf}{{\mathcal{F}}}
\newcommand{\scrg}{{\mathcal{G}}}
\newcommand{\scrh}{{\mathcal{H}}}
\newcommand{\scrk}{{\mathcal{K}}}
\newcommand{\scrl}{{\mathcal{L}}}
\newcommand{\scro}{{\mathcal{O}}}
\newcommand{\scrp}{{\mathcal{P}}}
\newcommand{\scrr}{{\mathcal{R}}}
\newcommand{\scrs}{{\mathcal{S}}}
\newcommand{\scrt}{{\mathcal{T}}}
\newcommand{\scrv}{{\mathcal{V}}}
\newcommand{\scrw}{{\mathcal{W}}}
\newcommand{\scrz}{{\mathcal{Z}}}

\renewcommand{\k}{{\mathbf{k}}}
\newcommand{\n}{{\mathbf{n}}}
\newcommand{\vv}{{\mathbf{v}}}
\newcommand{\bv}{{\mathbf{v}}}
\newcommand{\w}{{\mathbf{w}}}
\newcommand{\x}{{\mathbf{x}}}
\newcommand{\cc}{{\mathbf{c}}}
\newcommand{\zero}{{\mathbf{0}}}
\newcommand{\one}{{\mathbf{1}}}

\newcommand{\C}{{\mathbb C}}
\newcommand{\D}{{\mathbb D}}
\newcommand{\Dbar}{{\overline{\mathbb D}}}
\newcommand{\Z}{{\mathbb Z}}
\newcommand{\N}{{\mathbb N}}
\newcommand{\Q}{{\mathbb Q}}
\newcommand{\R}{{\mathbb R}}
\newcommand{\RR}{{\mathbb R}}

\newcommand{\ahat}{{\widehat{a}}}
\newcommand{\ctilde}{{\widetilde{c}}}
\newcommand{\cdoubletilde}{{\widetilde{\widetilde{c}}}}
\newcommand{\Zhat}{{\widehat{Z}}}
\newcommand{\balpha}{{\boldsymbol{\alpha}}}
\newcommand{\ginv}{{g^{\langle -1 \rangle}}}
\newcommand{\scrrtilde}{{\widetilde{\scrr}}}
\newcommand{\Rtilde}{{\widetilde{R}}}

\newcommand{\bX}{{\boldsymbol{X}}}
\newcommand{\bx}{{\boldsymbol{x}}}

\newcommand{\myle}{\preceq}
\newcommand{\myge}{\succeq}

%
% Hypergeometric functions
%
\newcommand{\ofo}{ {{}_1 \! F_1} }
\newcommand{\tfo}{ {{}_2 \! F_1} }
\newcommand{\ophio}{ {{}_1  \phi_1} }
\newcommand{\tphio}{ {{}_2  \phi_1} }

%
% Variants of \binom  (defined using the AMS "genfrac" command)
%
\newcommand{\stirlingsubset}[2]{\genfrac{\{}{\}}{0pt}{}{#1}{#2}}
\newcommand{\stirlingcycle}[2]{\genfrac{[}{]}{0pt}{}{#1}{#2}}
\newcommand{\assocstirlingsubset}[3]{{\genfrac{\{}{\}}{0pt}{}{#1}{#2}}_{\! \ge #3}}
\newcommand{\assocstirlingcycle}[3]{{\genfrac{[}{]}{0pt}{}{#1}{#2}}_{\ge #3}}
\newcommand{\genstirlingsubset}[4]{{\genfrac{\{}{\}}{0pt}{}{#1}{#2}}_{\! #3,#4}}
\newcommand{\euler}[2]{\genfrac{\langle}{\rangle}{0pt}{}{#1}{#2}}
\newcommand{\eulergen}[3]{{\genfrac{\langle}{\rangle}{0pt}{}{#1}{#2}}_{\! #3}}
\newcommand{\eulersecond}[2]{\left\langle\!\! \euler{#1}{#2} \!\!\right\rangle}
\newcommand{\eulersecondgen}[3]{{\left\langle\!\! \euler{#1}{#2} \!\!\right\rangle}_{\! #3}}
\newcommand{\binomvert}[2]{\genfrac{\vert}{\vert}{0pt}{}{#1}{#2}}
\newcommand{\qbinom}[3]{\genfrac{[}{]}{0pt}{}{#1}{#2}_{#3}}

% Array for subscripts

\newenvironment{sarray}{
             \textfont0=\scriptfont0
             \scriptfont0=\scriptscriptfont0
             \textfont1=\scriptfont1
             \scriptfont1=\scriptscriptfont1
             \textfont2=\scriptfont2
             \scriptfont2=\scriptscriptfont2
             \textfont3=\scriptfont3
             \scriptfont3=\scriptscriptfont3
           \renewcommand{\arraystretch}{0.7}
           \begin{array}{l}}{\end{array}}

\newenvironment{scarray}{
             \textfont0=\scriptfont0
             \scriptfont0=\scriptscriptfont0
             \textfont1=\scriptfont1
             \scriptfont1=\scriptscriptfont1
             \textfont2=\scriptfont2
             \scriptfont2=\scriptscriptfont2
             \textfont3=\scriptfont3
             \scriptfont3=\scriptscriptfont3
           \renewcommand{\arraystretch}{0.7}
           \begin{array}{c}}{\end{array}}

\clearpage

\section{Introduction}

Consider a formal power series of the form
\begin{equation}
   f(x,y) \;=\; \sum\limits_{n=0}^\infty \alpha_n \, x^n \, y^{n(n-1)/2}
 \label{def.f}
\end{equation}
where the coefficients $(\alpha_n)_{n=0}^\infty$ belong to a
commutative ring-with-identity-element $R$
and we impose the normalization $\alpha_0 = \alpha_1 = 1$.
We can regard $f$ as a formal power series in $y$
whose coefficients are polynomials in $x$, i.e.\ $f \in R[x][[y]]$.
Then, for any formal power series $X(y)$ with coefficients in $R$,
the composition $f(X(y),y)$ makes sense as a formal power series in $y$.
In particular, it is easy to see
--- either by the implicit function theorem for formal power series
\cite[p.~A.IV.37]{Bourbaki_90} \cite[Proposition~3.1]{Sokal_implicit}
or by a direct inductive argument ---
that there exists a unique formal power series $x_0(y) \in R[[y]]$
satisfying $f(x_0(y),y) = 0$,
which I~call the ``leading root'' of $f$.
Since $x_0(y)$ obviously has constant term $-1$,
it is convenient to write $x_0(y) = -\xi_0(y)$ where $\xi_0(y) = 1 + O(y)$.

Among the interesting series $f(x,y)$ of this type
are the ``partial theta function''
\cite[Chapter~13]{Andrews-Berndt_05} \cite[Chapter~6]{Andrews-Berndt_09}
\begin{equation}
   \Theta_0(x,y) \;=\; \sum\limits_{n=0}^\infty x^n \, y^{n(n-1)/2}
 \label{def.G}
\end{equation}
and the ``deformed exponential function''
\cite{Morris_72,Liu_98,Langley_00,Sokal_Fxy_asymptotic,%
Sokal_Fxy_conjectures,Sokal_Fxy_powerseries,Sokal_Fxy_hadamard}
\begin{equation}
   F(x,y) \;=\; \sum\limits_{n=0}^\infty {x^n \over n!} \, y^{n(n-1)/2}
   \;.
 \label{def.F}
\end{equation}
More generally one can consider the
rescaled three-variable Rogers--Ramanujan function \cite{Sokal_Fxy_powerseries}
\be
   \Rtilde(x,y,q)
   \;=\;
   \sum\limits_{n=0}^\infty
       {x^n \, y^{n(n-1)/2}  \over
        (1+q) (1+q+q^2) \,\cdots\, (1+q+\ldots+q^{n-1})}
   \;,
 \label{def.Rtilde}
\ee
which reduces to the foregoing when $q=0$ and $q=1$, respectively.

I have recently discovered empirically that the power series $\xi_0(y)$
has all nonnegative (in fact strictly positive) coefficients
in the first two cases, and more generally in the third case whenever $q > -1$.
More precisely, I have verified this for $\Theta_0$ and $F$
through orders $y^{6999}$ and $y^{899}$, respectively,
%% Fxy_general_x0_v2ab_Gxy_out_6500.m.README
%% Fxy_general_x0_v2a_Gxy.out
%% Fxy_general_x0_v2_Fxy.out
using a formula \cite{Sokal_Fxy_powerseries}
that relates $\xi_0(y)$ to the series expansion of $\log f(x,y)$.
For $\Rtilde$, I have proven \cite{Sokal_Fxy_powerseries} that $\xi_0(y,q)$ has
the form
\be
   \xi_0(y,q)  \;=\;  1 \,+\, \sum_{n=1}^\infty {P_n(q) \over Q_n(q)} \, y^n
 \label{eq.xi0.RR}
\ee
where
\be
   Q_n(q)
   \;=\;
   \prod_{k=2}^\infty (1+q+\ldots+q^{k-1})^{\lfloor n/{k \choose 2} \rfloor}
 \label{def.Qnq}
\ee
and $P_n(q)$ is a self-inversive polynomial in $q$ with integer coefficients;
and I have verified for $n \le 349$ that $P_n(q)$ has {\em two}\/
interesting positivity properties:
%% -rw-r--r--  1 as2a physics 113362 Nov 29 10:20 Fxy_general_x0_v2c_Rxy.out_f
%% -rw-r--r--  1 as2a physics 46686569 Nov 29 16:39 Fxy_general_x0_v2c_Rxy_out_350f_PPonly.m
\begin{itemize}
   \item[(a)] $P_n(q)$ has all nonnegative coefficients.
      Indeed, all the coefficients are strictly positive
      except $[q^1] \, P_5(q) = 0$.
   \item[(b)] $P_n(q) > 0$ for $q > -1$.
\end{itemize}
Of course, I conjecture that these properties hold for all $n$,
but I have (as yet) no proof.

The main purpose of this paper is to give a simple proof of the
coefficientwise positivity of $\xi_0(y)$
in the case of the partial theta function \reff{def.G}:

\begin{theorem}
   \label{thm1}
For the partial theta function \reff{def.G},
the formal power series
\be
   \quad
   \xi_0(y)  \;=\;
   1 + y + 2y^2 + 4y^3 + 9y^4 + 21y^5 + 52y^6 + 133y^7 + 351y^8 + 948y^9 +
      2610y^{10} +  \ldots \quad
\ee
has strictly positive coefficients.
\end{theorem}

In fact, with a bit more work one can prove a
pair of successively stronger results:

\begin{theorem}
   \label{thm2}
For the partial theta function \reff{def.G},
the formal power series
\be
   1/\xi_0(y)  \;=\;
   1 - y - y^2 - y^3 -2 y^4 - 4y^5 - 10y^6 - 25y^7 - 66y^8 - 178y^9 -
      490y^{10}  -  \ldots \quad
\ee
has strictly negative coefficients after the constant term 1.
\end{theorem}

\begin{theorem}
   \label{thm3}
For the partial theta function \reff{def.G},
the formal power series
\be
   1/\xi_0(y)^{2}  \;=\;
   1-2 y-y^2 \hphantom{-y^3} - y^4-2 y^5-7 y^6-18 y^7-50 y^8-138 y^9-386 y^{10}
   - \ldots
\ee
has strictly negative coefficients after the constant term 1
except for the vanishing coefficient of $y^3$.
\end{theorem}

\noindent
For further discussion of the relationship between these results,
see Section~\ref{sec.discussion}.

In addition, I have discovered empirically a vast strengthening of
Theorems~\ref{thm1} and \ref{thm2}.
Please note first that any power series
$g(y) = 1 + \sum_{n=1}^\infty a_n y^n \in \Z[[y]]$
can be written uniquely as an infinite product
$g(y) = \prod_{m=1}^\infty (1-y^m)^{-c_m}$
with coefficients $c_m \in \Z$.\footnote{
   See e.g.\ \cite[Theorem~10.3]{Andrews_86}.
   Some authors \cite{Bernstein_95,OEIS} \cite[pp.~20--21]{Sloane_95}
   call $(a_n)_{n=1}^\infty$ the {\em Euler transform}\/
   of $(c_m)_{m=1}^\infty$,
   and $(c_m)_{m=1}^\infty$ the {\em inverse Euler transform}\/
   of $(a_n)_{n=1}^\infty$.
   However, this should not be confused with an unrelated
   (and more widely used) ``Euler transformation'' of sequences,
   involving binomial coefficients.
}
We then have:

\begin{conjecture}
   \label{conj1.4}
For the partial theta function \reff{def.G},
when the formal power series $\xi_0(y)$ is written in the form
$\xi_0(y) = \prod\limits_{m=1}^\infty (1-y^m)^{-c_m}$,
% \be
%    \xi_0(y)  \;=\;  \prod_{m=1}^\infty (1-y^m)^{-c_m}
%    \;,
% \ee
the coefficient sequence
\be
   (c_m)_{m=1}^\infty
   \;=\;
   1, 1, 2, 4, 10, 23, 61, 157, 426, 1163, 3253, 9172, 26236, 75634, \ldots
\ee
is strictly positive ($c_m > 0$),
increasing ($\Delta c \ge 0$),
strictly convex ($\Delta^2 c > 0$),
and satisfies $\Delta^k c \ge 0$ for $k=3,4$.
[By contrast, the sequence $\Delta^5 c$ starts with $-3$.]
\end{conjecture}

\begin{conjecture}
   \label{conj1.5}
For the partial theta function \reff{def.G},
when the formal power series $\xi_0(y)$ is written in the form
$2 - 1/\xi_0(y) = \prod\limits_{m=1}^\infty (1-y^m)^{-c'_m}$,
% \be
%    2 \,-\, {1 \over \xi_0(y)}  \;=\;  \prod_{m=1}^\infty (1-y^m)^{-c'_m}
%    \;,
% \ee
the coefficient sequence
\be
   (c'_m)_{m=1}^\infty
   \;=\;
   1, 0, 0, 1, 2, 6, 15, 40, 110, 303, 853, 2419, 6950, 20110, \ldots
\ee
is nonnegative and convex.
\end{conjecture}

\noindent
I have verified these conjectures through order $y^{6999}$,
%% Use file Fxy_general_x0_v2ab_Gxy_out_7000_xi0only.m
%%   and program Fxy_general_Gxy_eulertrans.m
%% -rw-r--r-- 1 as2a physics 9658 Apr 12 14:22 Fxy_general_Gxy_eulertrans.out
%% -rw-r--r-- 1 as2a physics 2711 Apr 12 12:50 Fxy_general_Gxy_eulertrans.m
but I have no idea how to prove them.
Perhaps one should try to find a combinatorial interpretation
of the coefficients $(c_m)$ and $(c'_m)$.

The series $\xi_0(y)$ appears to possess one further striking property,
which I have again verified through order $y^{6999}$:

\begin{conjecture}
   \label{conj1.6}
For the partial theta function \reff{def.G},
the coefficient sequence of $\xi_0(y) = \sum_{n=0}^\infty a_n y^n$
is log convex, i.e.\ $a_{n-1} a_{n+1} \ge a_n^2$ for all $n \ge 1$.
\end{conjecture}

\noindent
A classic theorem of Kaluza \cite{Kaluza_28}
relates Conjecture~\ref{conj1.6} to Theorem~\ref{thm2}:
namely, if the coefficient sequence $(a_n)_{n=0}^\infty$
of a formal power series $f$ is strictly positive and log convex,
then $1/f$ has nonpositive coefficients after the constant term;
and if in addition $a_0 a_2 > a_1^2$, then $1/f$ has
strictly negative coefficients after the constant term.\footnote{
   The assertion about strict negativity
   is not explicitly stated by Kaluza \cite{Kaluza_28},
   but it follows easily from his proof.
   See also \cite[Lemma~2.2]{Krattenthaler_11}.
}
But it is easily seen that the converse does not hold.\footnote{
   For instance, let $f(y) = 1/(1-y-cy^2)$;
   then $1/f$ has nonpositive coefficients after the constant term
   whenever $c \ge 0$;
   but the coefficients of $f$ are log convex only when $c=0$.
}
So Conjecture~\ref{conj1.6}, if true,
is a strengthening of Theorem~\ref{thm2}.

\bigskip

The plan of this paper is as follows:
I begin (Section~\ref{sec.identities})
by recalling two identities for the partial theta function,
which will play a central role in the proofs of
Theorems~\ref{thm1}--\ref{thm3}.
%% I then give the proofs of Theorem~\ref{thm1} (Section~\ref{sec.proof})
%% and Theorem~\ref{thm2} (Section~\ref{sec.proof2}).
I then give, in successive sections,
the proofs of Theorem~\ref{thm1}--\ref{thm3}
(Sections~\ref{sec.proof}--\ref{sec.proof3}).
% Next I state and prove two generalizations (Section~\ref{sec.generalization}).
% Then
Next I state and prove some identities for the
three-variable Rogers--Ramanujan function \reff{def.Rtilde}
that may turn out to be useful in proving the conjectures
concerning its leading root (Section~\ref{sec.Rxyq}).
Finally, I place Theorems~\ref{thm1}--\ref{thm3}
in a more general context \cite{Scott-Sokal_powerseries}
and mention some stronger properties
possessed by the power series $\xi_0(y)$
for the cases \reff{def.G}--\reff{def.Rtilde}
that appear empirically to be true (Section~\ref{sec.discussion}).

A {\sc Mathematica} file {\tt partialtheta\_xi0.m}
containing the series $\xi_0(y)$ for the partial theta function
through order $y^{6999}$ is available as an ancillary file
with the preprint version of this paper at \url{arXiv.org}.

\section{Identities for the partial theta function}   \label{sec.identities}

In this section we recall a pair of identities
for the partial theta function \reff{def.G}
that will serve as the foundation for our proofs of
Theorems~\ref{thm1}--\ref{thm3}.
We use the standard notation $(a;q)_n = \prod_{j=0}^{n-1} (1-aq^j)$
and $(a;q)_\infty = \prod_{j=0}^{\infty} (1-aq^j)$.

\begin{lemma}
   \label{lemma.identities}
The partial theta function \reff{def.G} satisfies
\begin{eqnarray}
   %% \sum\limits_{n=0}^\infty  x^n \, y^{n(n-1)/2}
   \Theta_0(x,y)
   & = &
   (y;y)_\infty \, (-x;y)_\infty \,
      \sum_{n=0}^\infty  {y^n  \over  (y;y)_n \, (-x;y)_n}
    \label{eq.lemma.identity1}  \\[2mm]
   %% \sum\limits_{n=0}^\infty  x^n \, y^{n(n-1)/2}
   \Theta_0(x,y)
   & = &
   (-x;y)_\infty \,
      \sum_{n=0}^\infty  {(-x)^n \, y^{n^2}  \over  (y;y)_n \, (-x;y)_n}
    \label{eq.lemma.identity2}
\end{eqnarray}
as formal power series
and as analytic functions on $(x,y) \in \C \times \D$.\footnote{
   Here $\D$ denotes the open unit disc in $\C$.
   The right-hand sides of
   \reff{eq.lemma.identity1} and \reff{eq.lemma.identity2}
   have removable singularities at $x=-y^{-k}$ ($k=0,1,2,\ldots$).
   To see that these singularities are indeed removable,
   just rewrite $(-x;y)_\infty / (-x;y)_n$ as $(-xy^n;y)_\infty$.
}
\end{lemma}

In order to make this paper self-contained for readers who (like myself!)\ 
are not experts in $q$-series,
we provide here an easy proof of \reff{eq.lemma.identity1}
that uses nothing more than Euler's first and second identities
%% \cite[Corollary~2.2]{Andrews_76}
\cite[eqs.~(1.3.15) and (1.3.16)]{Gasper_04}
\begin{eqnarray}
   {1 \over (t;q)_\infty}
   & = &
   \sum\limits_{n=0}^\infty  {t^n \over (q;q)_n}
      \label{eq.euler1}    \\[2mm]
   (t;q)_\infty
   & = &
   \sum\limits_{n=0}^\infty  {(-t)^n \, q^{n(n-1)/2} \over (q;q)_n}
      \label{eq.euler2}
\end{eqnarray}
valid for $(t,q) \in \D \times \D$ and $(t,q) \in \C \times \D$, respectively.

%%\proofof{\reff{eq.lemma.identity1}}
\bigskip\noindent
{\sc Proof of \reff{eq.lemma.identity1}} \cite{Chen_09b,Alladi_09b}.\ \ 
Write
\be
   %% \sum\limits_{n=0}^\infty  x^n \, y^{n(n-1)/2}
   \Theta_0(x,y)
   \;=\;
   \sum\limits_{n=0}^\infty x^n \, y^{n(n-1)/2}
         \: {(y;y)_\infty  \over (y;y)_n \, (y^{n+1};y)_\infty}
 \label{eq.proof.eq.lemma.identity1.1}
\ee
and insert Euler's first identity for $1/(y^{n+1};y)_\infty$:
we obtain
\begin{subeqnarray}
   \Theta_0(x,y)
   & = &
   (y;y)_\infty
   \sum\limits_{n=0}^\infty  {x^n \, y^{n(n-1)/2}  \over  (y;y)_n}
   \sum\limits_{k=0}^\infty  {y^{(n+1)k}  \over (y;y)_k}
             \qquad \\[1mm]
   & = &
   (y;y)_\infty
   \sum\limits_{k=0}^\infty
       {y^k  \over (y;y)_k}
   \sum\limits_{n=0}^\infty  {(xy^k)^n \, y^{n(n-1)/2}  \over  (y;y)_n}
             \qquad \\[1mm]
   & = &
   (y;y)_\infty
   \sum\limits_{k=0}^\infty
       {y^k  \over (y;y)_k}
       \: (-xy^k; y)_\infty  \quad\hbox{by Euler's second identity}
             \qquad \\[1mm]
   & = &
   (y;y)_\infty \, (-x;y)_\infty \,
   \sum\limits_{k=0}^\infty
       {y^k  \over (y;y)_k \, (-x; y)_k}
   \;.
 \label{eq.proof.eq.lemma.identity1.2}
\end{subeqnarray}
\qed

And here is an easy proof of
both \reff{eq.lemma.identity1} and \reff{eq.lemma.identity2}
that uses only Heine's first and second transformations
%% \cite[Corollary~2.3 and eq.~(3.3.13)]{Andrews_76}
\cite[eqs.~(1.4.1) and (1.4.5)]{Gasper_04}
\begin{eqnarray}
   \tphio(a,b;c;q,z)
   & = &
   {(b;q)_\infty \, (az;q)_\infty \over (c;q)_\infty \, (z;q)_\infty}
   \;
   \tphio(c/b,z;az;q,b)
        \label{eq.heine1} \\[2mm]
   \tphio(a,b;c;q,z)
   & = &
   {(c/a;q)_\infty \, (az;q)_\infty \over (c;q)_\infty \, (z;q)_\infty}
   \;
   \tphio(abz/c,a;az;q,c/a)
        \label{eq.heine2}
\end{eqnarray}
for the basic hypergeometric function
\be
   \tphio(a,b;c;q,z)
   \;=\;
   \sum_{n=0}^\infty {(a;q)_n \, (b;q)_n \over (q;q)_n \, (c;q)_n} \: z^n
   \;.
 \label{def.tphio}
\ee
Here \reff{eq.heine1} is valid when $|q| < 1$, $|z| < 1$ and $0 < |b| < 1$,
while \reff{eq.heine2} is valid when $|q| < 1$, $|z| < 1$ and $0 < |c| < |a|$.

\bigskip\noindent
{\sc Proof of \reff{eq.lemma.identity1} and \reff{eq.lemma.identity2}}
 \cite{Andrews_07}.\ \ 
In \reff{eq.heine1} and \reff{eq.heine2}, set $b=q$ and $z=-x/a$,
then take $a \to\infty$ and $c \to 0$;
we obtain \reff{eq.lemma.identity1} and \reff{eq.lemma.identity2}
with $y$ renamed as $q$.
\qed

\bigskip

{\bf Remarks.}
Identity \reff{eq.lemma.identity1}
goes back to Heine in 1847 \cite[bottom p.~306]{Heine_1847},
who derived it (as here) as a limiting case of his
fundamental transformation \reff{eq.heine1}.\footnote{
   Heine makes the change of variables $x=-zq$ and $y=q^2$.
   The formula in \cite[bottom p.~306]{Heine_1847}
   has a typographical error in which the factor $y^n$ (= $q^{2n}$)
   in the numerator of the right-hand side is inadvertently omitted.
   The correct formula can be found in the 1878 edition
   of Heine's book \cite[p.~107]{Heine_1878}.
}
In the modern literature it can be found in Fine \cite[eq.~(7.32)]{Fine_88}.
%% {\bf Where else?????  Surely there are references in-between
%%    1847 and 1988!!!!}

I don't know who first found identity \reff{eq.lemma.identity2};
I would be grateful to any reader who can supply a reference.
I first learned \reff{eq.lemma.identity2}
from the paper of Andrews and Warnaar \cite[eq.~(2.1)]{Andrews_07},
but it is surely much older.
%% but a more general identity can be found in the earlier papers
%% of Bhargava and Adiga \cite{Bhargava_86} and Srivastava \cite{Srivastava_87}.

The elementary proof of \reff{eq.lemma.identity1} given here
is in essence that given recently by
Chen and Xia \cite[eq.~(2.10)]{Chen_09b}
and Alladi \cite[p.~second proof of~(1.6)]{Alladi_09b}.\footnote{
   Alladi's eq.~(1.6) is equivalent to our \reff{eq.lemma.identity1}
   under the substitutions $x=-aq$ and $y=q^2$.
}
Our proof of \reff{eq.lemma.identity1} and \reff{eq.lemma.identity2}
using Heine's transformations follows
Andrews and Warnaar \cite[eq.~(2.1)]{Andrews_07}\footnote{
   See also Andrews \cite[proof of Theorem~1]{Andrews_05b}
   for this proof of \reff{eq.lemma.identity1}.
},
but at least for \reff{eq.lemma.identity1}
the argument goes back to Heine himself \cite[p.~306]{Heine_1847}.
Note also that if one takes this latter proof of \reff{eq.lemma.identity1}
and inserts in it the standard {\em proof}\/ of Heine's first transformation
\cite[sec.~1.4]{Gasper_04},
one obtains the elementary proof of \reff{eq.lemma.identity1}.
%% given in
%% \reff{eq.proof.eq.lemma.identity1.1}/\reff{eq.proof.eq.lemma.identity1.2}
%%   above.

A combinatorial proof of \reff{eq.lemma.identity1}
was given recently by Yee \cite[Theorem~2.1]{Yee_10},
and combinatorial proofs of both \reff{eq.lemma.identity1}
and the equality \reff{eq.lemma.identity1}=\reff{eq.lemma.identity2}
were given recently by Kim \cite[Section~2]{Kim_10}.

Many generalizations of \reff{eq.lemma.identity1}/\reff{eq.lemma.identity2},
with additional parameters, are known.
%% One such generalization is given in Lemma~\ref{lemma.identities.ophio} below.
For instance, \reff{eq.lemma.identity1}/\reff{eq.lemma.identity2}
can be extended from the partial theta function to more general
basic hypergeometric functions $\ophio$.\footnote{
   For the case of \reff{eq.lemma.identity2},
   this generalization can be found in papers of
   Bhargava and Adiga \cite{Bhargava_86} and Srivastava \cite{Srivastava_87}.
   A special case of this generalization
   can be found in Ramanujan's second notebook
   \cite[Entry~9 in Chapter~16]{Ramanujan_57} \cite[p.~18]{Berndt_91}
   and again in a page published with the lost notebook
   \cite[p.~362]{Ramanujan_lost} \cite[Entry~1.6.1]{Andrews-Berndt_09}.
   A combinatorial proof of this special case was recently given by
   Berndt, Kim and Yee \cite[Theorem~5.6]{Berndt_10}.
}
Another generalization of \reff{eq.lemma.identity1}
appears in Ramanujan's lost notebook
\cite[p.~40]{Ramanujan_lost} \cite[Entry~6.3.1]{Andrews-Berndt_09};
it was proven by Andrews \cite[Section~4]{Andrews_79}
and recently re-proven combinatorially by Kim \cite[Section~4]{Kim_10}.
An even more general formula was proven
subsequently by Andrews \cite[Section~3]{Andrews_81},
with a later simplification
and further generalization by R.P.~Agarwal \cite{Agarwal_84};
see also \cite[Sections~6.2 and 6.3]{Andrews-Berndt_09}.
A formula generalizing the equality
\reff{eq.lemma.identity1}=\reff{eq.lemma.identity2}
appears in Ramanujan's lost notebook
\cite[p.~40]{Ramanujan_lost} \cite[Entry~1.6.7]{Andrews-Berndt_09}
and has an easy $q$-series proof \cite[p.~27]{Andrews-Berndt_09};
a combinatorial proof was recently given by Kim \cite[Section~4]{Kim_10}.

A very beautiful formula for the sum of two partial theta functions,
which generalizes both \reff{eq.lemma.identity1}
and the Jacobi triple product identity,
was found by Warnaar \cite[Theorem~1.5]{Warnaar_03}.
A closely related identity for the product of two partial theta functions,
which also generalizes \reff{eq.lemma.identity1},
was found by Andrews and Warnaar \cite[Theorem~1.1]{Andrews_07}
and recently re-proven combinatorially by Kim \cite[Section~3]{Kim_10};
see also \cite[Section~6.6]{Andrews-Berndt_09}.
%% {\bf Take a closer look also at
%%    Andrews \cite[p.~29, Examples~5 and 6]{Andrews_76}
%%    and the Andrews 1966a,c references cited on p.~31.}

Finally, Andrews \cite[Theorem~5]{Andrews_09b}
has recently proven a finite-sum generalization of \reff{eq.lemma.identity1}:
\be
   \sum_{n=0}^N {x^n \, y^{n(n-1)/2} \over (y;y)_{N-n}}
   \;=\;
   (-x;y)_N \, \sum_{n=0}^N  {y^n  \over  (y;y)_n \, (-x;y)_n}
   \;.
 \label{eq.finite.identity1}
\ee
Likewise, by using \cite[Corollary~3]{Andrews_09b}
with $\alpha=q$, $\tau=-x/\beta$ and taking
$\beta \to \infty$ and $\gamma \to 0$,
one can derive a finite-sum generalization of \reff{eq.lemma.identity2}:
\be
   \sum_{n=0}^N {x^n \, y^{n(n-1)/2} \over (y;y)_{N-n}}
   \;=\;
   (-x;y)_N \, \sum_{n=0}^N
       {(-x)^n \, y^{n^2}  \over  (y;y)_n \, (-x;y)_n \, (y;y)_{N-n}}
   \;.
 \label{eq.finite.identity2}
\ee
See also \cite{Rowell_09} for a combinatorial proof of the
finite Heine transformation that underlies
\reff{eq.finite.identity1} and \reff{eq.finite.identity2}.

\section{Proof of Theorem~\ref{thm1}}   \label{sec.proof}

The proof of Theorem~\ref{thm1} can be based on
either \reff{eq.lemma.identity1} or \reff{eq.lemma.identity2}.
For concreteness let us use \reff{eq.lemma.identity1},
which we rewrite as
\be
   \Theta_0(x,y)
   \;=\;
   %% \sum_{n=0}^\infty y^n \, (y^{n+1};y)_\infty \, (-xy^n;y)_\infty  \\
   %% & = &
   %% (y;y)_\infty \, (-x;y)_\infty \,
   %%   \sum_{n=0}^\infty  {y^n  \over  (y;y)_n \, (-x;y)_n}   \\
   %% & = &
   (y;y)_\infty \, (-xy;y)_\infty
      \left[ 1 + x + \sum_{n=1}^\infty  {y^n  \over  (y;y)_n \, (-xy;y)_{n-1}}
      \right]  .
 \label{eq.identity.c}
\ee
So $\Theta_0(-\xi_0(y),y) = 0$ is equivalent to
\begin{subeqnarray}
   \xi_0(y)
   & = &
   1 \,+\, \sum_{n=1}^\infty  {y^n  \over  (y;y)_n \, (y\xi_0(y);y)_{n-1}}  \\
   & = &
   1 \,+\, \sum_{n=1}^\infty
       {y^n  \over  \prod\limits_{j=1}^n (1-y^j)
                    \prod\limits_{j=1}^{n-1} [1-y^j \xi_0(y)] }
   \;\,.
 \label{eq.identity2}
\end{subeqnarray}
This formula can be used iteratively to determine $\xi_0(y)$,
and in particular to prove the strict positivity of its coefficients:

\begin{proposition}
   \label{prop2}
Define the map $\scrf \colon\, \Z[[y]] \to \Z[[y]]$ by
\be
   (\scrf \xi)(y)
   \;=\;
   1 \,+\, \sum_{n=1}^\infty
       {y^n  \over  \prod\limits_{j=1}^n (1-y^j)
                    \prod\limits_{j=1}^{n-1} [1-y^j \xi(y)] }
   \;\,,
 \label{def.scrf}
\ee
and define a sequence $\xi_0^{(0)}, \xi_0^{(1)}, \ldots \in \Z[[y]]$ by
$\xi_0^{(0)} = 1$ and $\xi_0^{(k+1)} = \scrf \xi_0^{(k)}$.
Then
\be
   \xi_0^{(0)} \,\preceq\, \xi_0^{(1)} \,\preceq\, \xi_0^{(2)}
      \,\preceq\, \ldots \,\preceq\, \xi_0
 \label{eq.xi0xi1etc}
\ee
(where $f \preceq g$ denotes $[y^n] f(y) \le [y^n] g(y)$ for all $n$)
and
\be
   \xi_0^{(k)}(y) \;=\; \xi_0(y) \,+\, O(y^{3k+1})
   \;.
\ee
In particular, $\lim\limits_{k\to\infty} \xi_0^{(k)}(y) = \xi_0(y)$
in the sense of convergence of formal power series
(i.e.\ every coefficient eventually stabilizes at its limit),
and $\xi_0(y)$ has strictly positive coefficients.
\end{proposition}

\proof
If $f(y)$ and $g(y)$ are formal power series
satisfying $0 \preceq f \preceq g$,
then it is easy to see that
$\prod\limits_{j=1}^{n-1} [1-y^j f(y)]^{-1}
 \preceq
 \prod\limits_{j=1}^{n-1} [1-y^j g(y)]^{-1}$
and hence $ 0 \preceq \scrf f \preceq \scrf g$.
Applying this repeatedly to the obvious inequality
$0 \preceq \xi_0^{(0)} \preceq \xi_0^{(1)}$,
we obtain
$\xi_0^{(0)} \preceq \xi_0^{(1)} \preceq \xi_0^{(2)} \preceq \ldots\;$.

Likewise, if $f(y)$ and $g(y)$ are formal power series
satisfying $f(y) - g(y) = O(y^\ell)$ for some $\ell \ge 0$,
then it is not hard to see that
$(\scrf f)(y) - (\scrf g)(y) = O(y^{\ell+3})$
[coming from the $n=2$ term in \reff{def.scrf}
 and the $j=1$ factor in the second product].
Applying this repeatedly to the obvious fact
$\xi_0^{(1)}(y) - \xi_0^{(0)}(y) = O(y)$,
we obtain
$\xi_0^{(k+1)}(y) - \xi_0^{(k)}(y) = O(y^{3k+1})$.
It follows that $\xi_0^{(k)}(y)$ converges as $k\to\infty$
(in the topology of formal power series)
to a limiting series $\xi_0^{(\infty)}(y)$,
and that this limiting series satisfies
$\scrf \xi_0^{(\infty)} = \xi_0^{(\infty)}$.
But this means, by \reff{eq.identity.c}/\reff{eq.identity2},
that $\xi_0^{(\infty)}(y) = \xi_0(y)$.
It also follows that $\xi_0^{(k)}(y) = \xi_0(y) + O(y^{3k+1})$.

Since $\xi_0^{(1)}(y)$ manifestly has strictly positive coefficients,
it follows from \reff{eq.xi0xi1etc} that $\xi_0(y)$ also has
strictly positive coefficients.
\qed

{\bf Remarks.}
1.
By a slightly more refined version of the same argument,
one can prove inductively that
\be
   \xi_0^{(k+1)}(y) - \xi_0^{(k)}(y)
   \;=\;
   y^{3k+1} \,+\, (4k+2) y^{3k+2} \,+\, (4k+1)(2k+3) y^{3k+3} \,+\, O(y^{3k+4})
\ee
for $k \ge 1$,
and hence that
\be
   \xi_0(y) - \xi_0^{(k)}(y)
   \;=\;
   y^{3k+1} \,+\, (4k+2) y^{3k+2} \,+\, (4k+1)(2k+3) y^{3k+3} \,+\, O(y^{3k+4})
\ee
for $k \ge 1$.
%% See my computations 4/29/11 at ENS-Paris, which I have scanned.
%% One has to consider contributions from n=2,3,4 in the series for \scrf.

2. The series
\begin{subeqnarray}
   \xi_0^{(1)}(y)
   & = &
   1 \,+\, \sum_{n=1}^\infty  {y^n  \over  (y;y)_n \, (y;y)_{n-1}}
   \;=\;
   1 \,+\, {1 - \Theta_0(-y,y)   \over (y;y)_\infty^2}
         \\[2mm]
   & = &
   1+y+2 y^2+4 y^3+8 y^4+15 y^5+27 y^6+47 y^7+79 y^8 + \ldots \quad
 \slabel{eq.stacks.b}
\end{subeqnarray}
enumerates weakly unimodal sequences of positive integers
(also called ``stacks'' or ``stack polyominoes'') by total weight
\cite{Auluck_51,Wright_71} \cite[Section~2.5]{Stanley_86}
\cite[sequence A001523]{OEIS}.
It would be interesting to seek combinatorial interpretations
of $\xi_0^{(k)}(y)$ for $k \ge 2$, or at least of $\xi_0(y)$.\footnote{
   {\bf Note added:}
   Thomas Prellberg \cite{Prellberg_11} has recently found
   a combinatorial interpretation of $\xi_0(y)$ and $\xi_0^{(k)}(y)$
   in terms of rooted trees enriched by stack polyominoes,
   using results from \cite{Prellberg_95} and \cite[Chapter~3]{Bergeron_98}.
      \label{note_Prellberg_1}
}

3. Empirically I have observed that the $\xi_0^{(k)}$
obey inequalities stronger than \reff{eq.xi0xi1etc},
namely $\xi_0^{(k)} / \xi_0^{(k-1)} \myge 1$ for $k \ge 1$.
I have verified this through order $y^{500}$ for $1 \le k \le 20$,
%% -rw-r--r-- 1 as2a physics 6749 Jan 14 07:48 leadingroot_partialtheta_fflist.out
%% -rw-r--r-- 1 as2a physics 5419167 Jan 14 07:47 leadingroot_partialtheta_fflist_out.m
%% -rw-r--r-- 1 as2a physics    3012 Jan  8 00:43 leadingroot_partialtheta_fflist.m
but I do not see how to prove it.
If true, this exhibits $\xi_0(y)$ as an infinite product of nonnegative series
$\xi_0^{(k)}(y) / \xi_0^{(k-1)}(y)$,
reminiscent of but different from Conjecture~\ref{conj1.4}.

4. The recursion $\xi_0^{(k+1)} = \scrf \xi_0^{(k)}$
could alternatively have been started with
$\xi_0^{(0)} = 0$ instead of $\xi_0^{(0)} = 1$.
The only difference is that we would then have
$\xi_0^{(k)}(y) - \xi_0(y) = O(y^{3k})$ instead of $O(y^{3k+1})$.
In this case
\begin{subeqnarray}
   \xi_0^{(1)}(y)
   & = &
   \sum_{n=0}^\infty  {y^n  \over  (y;y)_n}
   \;=\;
   {1 \over (y;y)_\infty}
   \;=\;
   \sum_{n=0}^\infty p(n) \, y^n
           \\[2mm]
   & = &
   1+y+2 y^2+3 y^3+5 y^4+7 y^5+11 y^6+15 y^7+22 y^8+ \ldots \quad
\end{subeqnarray}
is the generating function for all partitions of the integer $n$.
%% \cite[sequence A000041]{OEIS}.
Perhaps $\xi_0^{(k)}(y)$ for $k \ge 2$ have a simpler interpretation
with this choice of $\xi_0^{(0)}$.\footnote{
   {\bf Note added:}
   Thomas Prellberg \cite{Prellberg_11} has found
   a combinatorial interpretation of $\xi_0^{(k)}(y)$
   also for this choice of $\xi_0^{(0)}$.
      \label{note_Prellberg_2}
}

Furthermore, with this choice of $\xi_0^{(0)}$
we have empirically not only $\xi_0^{(k)} / \xi_0^{(k-1)} \myge 1$
for $k \ge 2$, but in fact
$\xi_0^{(k)}(y) / \xi_0^{(k-1)}(y) = \prod_{m=1}^\infty (1-y^m)^{-c_m^{(k)}}$
with nonnegative coefficients $c_m^{(k)}$.
I have verified this through order $y^{500}$ for $2 \le k \le 20$.
%% -rw-r--r-- 1 as2a physics 6749 Jan 14 07:48 leadingroot_partialtheta_fflist.out
%% -rw-r--r-- 1 as2a physics 5419167 Jan 14 07:47 leadingroot_partialtheta_fflist_out.m
%% -rw-r--r-- 1 as2a physics    3012 Jan  8 00:43 leadingroot_partialtheta_fflist.m
If true, this implies Conjecture~\ref{conj1.4}.

5.  If we use \reff{eq.lemma.identity2} instead of \reff{eq.lemma.identity1},
then we are led to the recursion based on the map
$\scrg \colon\, \Z[[y]] \to \Z[[y]]$ defined by
\be
   (\scrg \xi)(y)
   \;=\;
   1 \,+\, \sum_{n=1}^\infty
       {\xi(y)^n \, y^{n^2}  \over  \prod\limits_{j=1}^n (1-y^j)
                    \prod\limits_{j=1}^{n-1} [1-y^j \xi(y)] }
   \;\,.
 \label{def.scrg}
\ee
Using $\xi_0^{(0)} = 1$, we have for this map the slower convergence
$\xi_0^{(k)}(y) - \xi_0(y) = O(y^{k})$
[coming from the $\xi(y)^n$ factor in the numerator of the $n=1$ term
 in \reff{def.scrg}].
In this case the series
\begin{subeqnarray}
   \xi_0^{(1)}(y)
   & = &
   1 \,+\, \sum_{n=1}^\infty  {y^{n^2}  \over  (y;y)_n \, (y;y)_{n-1}}
   \;=\;
   1 \,+\, {1 - \Theta_0(-y,y)   \over (y;y)_\infty}
       \\[2mm]
   & = &
   1+y+y^2+y^3+2 y^4+3 y^5+5 y^6+7 y^7+10 y^8+ \ldots \quad
\end{subeqnarray}
enumerates $n$-stacks with strictly receding walls
\cite{Auluck_51,Wright_71} \cite[sequence A001522]{OEIS}.
Once again we have empirically
$\xi_0^{(k)} / \xi_0^{(k-1)} \myge 1$ for $k \ge 1$;
I have verified this through order $y^{2000}$ for $1 \le k \le 20$.
%% -rw-r--r-- 1 as2a physics     13529 Jan 12 00:07 leadingroot_partialtheta_gglist.out
%% -rw-r--r-- 1 as2a physics 134801305 Jan 11 23:16 leadingroot_partialtheta_gglist_out.m
%% -rw-r--r-- 1 as2a physics      2853 Jan  8 00:39 leadingroot_partialtheta_gglist.m
Furthermore, for this map taking $\xi_0^{(0)} = 0$ yields $\xi_0^{(1)} = 1$,
so we obtain the {\em same}\/ sequence (shifted by one)
with both initial conditions.\footnote{
   {\bf Note added:}
   Thomas Prellberg \cite{Prellberg_11} has found
   a combinatorial interpretation also for these $\xi_0^{(k)}(y)$.
      \label{note_Prellberg_3}
}
\qed

\medskip

It is useful to abstract what we have done here
(see \cite{Sokal_Fxy_powerseries} for details and extensions).
Consider a formal power series
(with coefficients in a commutative ring-with-identity-element $R$)
\be
   f(x,y) \;=\;  \sum\limits_{n=0}^\infty a_n(y) \, x^n
 \label{eq.f.general}
\ee
where
\begin{itemize}
   \item[(a)] $a_0(0) = a_1(0) = 1$;
   \item[(b)] $a_n(0) = 0$ for $n \ge 2$; and
   \item[(c)] $a_n(y) = O(y^{\nu_n})$ with
      $\lim\limits_{n\to\infty} \nu_n = \infty$.
\end{itemize}
Then it is easy to see that there exists
a unique formal power series $\xi_0(y)$ with coefficients in $R$
satisfying $f(-\xi_0(y),y) = 0$, and it has constant term 1.
Let us rearrange $f(-\xi_0(y),y) = 0$ as
\be
   \xi_0(y)
   \;=\;
   1 \,+\, \sum_{n=0}^\infty (-1)^n \, \widehat{a}_n(y) \, \xi_0(y)^n
   \;,
 \label{eq.recursion.xi0}
\ee
where $\widehat{a}_n(y)$ is defined by
\be
   \widehat{a}_n(y)
   \;=\;
   \begin{cases}
      a_n(y) - 1    & \text{\rm for $n=0,1$} \\
      a_n(y)        & \text{\rm for $n \ge 2$}
   \end{cases}
 \label{def.ahat}
\ee
Now suppose that the ring $R$ carries a partial order
compatible with the ring structure
(typically we will have $R=\R$, $\Q$ or $\Z$)
and that
\be
   (-1)^n \widehat{a}_n(y)  \;\myge\; 0  \qquad\hbox{for all } n \ge 0
   \;,
 \label{eq.prop.simple.scrs.hypothesis}
\ee
where $f(y) \myge 0$ means that $f$ has all nonnegative coefficients.
Then the recursion argument used in Proposition~\ref{prop2},
applied to \reff{eq.recursion.xi0},
shows that
$\xi_0(y) \myge 1 + \sum_{n=0}^\infty (-1)^n \, \widehat{a}_n(y)$.
The case treated here was
\be
   f(x,y)
   \;=\;
   {\Theta_0(x,y) \over (y;y)_\infty \, (-xy;y)_\infty}
   \;=\;
   1 + x + \sum_{n=1}^\infty  {y^n  \over  (y;y)_n \, (-xy;y)_{n-1}}
   \;.
\ee
The value of the identity \reff{eq.lemma.identity1}
or \reff{eq.lemma.identity2} for our purposes
is that powers of $x$ on the left-hand side
are transformed into powers of $-x$ on the right-hand side,
so that \reff{eq.prop.simple.scrs.hypothesis} holds for the latter.

\section{Proof of Theorem~\ref{thm2}}   \label{sec.proof2}

In this section we prove Theorem~\ref{thm2}
on the strict negativity of the coefficients of $\xi_0(y)^{-1}$
after the constant term 1.
It is convenient to state and prove first an abstract result of this form
\cite{Sokal_Fxy_powerseries};
then we verify the hypotheses of this abstract result in our specific case.

\begin{proposition}
  \label{prop.scrs-1.abstract}
Consider a formal power series
(with coefficients in a partially ordered commutative ring $R$)
\be
   f(x,y) \;=\;  \sum\limits_{n=0}^\infty a_n(y) \, x^n
\ee
where
\begin{itemize}
   \item[(a)] $a_0(0) = a_1(0) = 1$;
   \item[(b)] $a_n(0) = 0$ for $n \ge 2$; and
   \item[(c)] $a_n(y) = O(y^{\nu_n})$ with
      $\lim\limits_{n\to\infty} \nu_n = \infty$.
\end{itemize}
Let $\xi_0(y)$ be the unique power series satisfying $f(-\xi_0(y),y) = 0$.
Suppose that
\be
   1 \,-\, {a_1(y) \over a_0(y)}  \;\myge\;  0
\ee
and that
\be
   (-1)^n \, {a_n(y) \over a_0(y)}  \;\myge\;  0
   \qquad\hbox{for all } n \ge 2 \;.
\ee
Then
\be
   \xi_0(y)^{-1}
   \;\myle\;
   {a_1(y) \over a_0(y)}  \,-\,
      \sum_{n=2}^\infty (-1)^n \, {a_n(y) \over a_0(y)}
   \;\myle\;
   1
   \;.
 \label{eq.cor.simple.scrs.b}
\ee
\end{proposition}

\proof
Start from the equation
$\sum\limits_{n=0}^\infty (-1)^n a_n(y) \, \xi_0(y)^n = 0$,
divide by $a_0(y) \xi_0(y)$, and bring $\xi_0(y)^{-1}$ to the left-hand side:
we have
\be
   \xi_0(y)^{-1}
   \;=\;
   {a_1(y) \over a_0(y)}  \,-\,
      \sum_{n=2}^\infty (-1)^n \, {a_n(y) \over a_0(y)} \, \xi_0(y)^{n-1}
   \;.
 \label{eq.xi0.recursion.b}
\ee
Now write $\xi_0(y)^{-1} = 1 - \psi(y)$: we obtain
\be
   \psi(y)
   \;=\;
   1 \,-\, {a_1(y) \over a_0(y)}  \,+\,
      \sum_{n=2}^\infty (-1)^n \, {a_n(y) \over a_0(y)} \, [1-\psi(y)]^{-(n-1)}
   \;.
 \label{eq.xi0.recursion.b.bis}
\ee
By hypothesis \reff{eq.xi0.recursion.b.bis} is of the form
\be
   \psi(y)
   \;=\;
   b_1(y) \,+\, \sum_{n=2}^\infty b_n(y) \, [1-\psi(y)]^{-(n-1)}
\ee
where $b_n(y) \myge 0$ and $b_n(y) = O(y)$ for all $n \ge 1$.
An iterative argument as in the proof of Proposition~\ref{prop2}
then proves that $\psi(y) \myge 0$ and in fact
\be
   \psi(y)
   \;\myge\;
   1 \,-\, {a_1(y) \over a_0(y)}  \,+\,
      \sum_{n=2}^\infty (-1)^n \, {a_n(y) \over a_0(y)}
   \;.
\ee
\qed

\proofof{Theorem~\ref{thm2}}
This time we find it convenient to use
\reff{eq.lemma.identity2} instead of \reff{eq.lemma.identity1}.
We therefore apply Proposition~\ref{prop.scrs-1.abstract} to the power series
\begin{subeqnarray}
   f(x,y)
   \;=\;
   {\Theta_0(x,y)  \over  (-xy;y)_\infty}
   & = &
   1 + x + \sum_{n=1}^\infty  {(-x)^n \, y^{n^2}  \over
                               (y;y)_n \, (-xy;y)_{n-1}}
          \\[2mm]
   & = &
   1 + x - {xy \over 1-y} +
       \sum_{n=2}^\infty  {(-x)^n \, y^{n^2}  \over (y;y)_n \, (-xy;y)_{n-1}}
   \;. \qquad
 \slabel{eq.proof.thm2.b}
\end{subeqnarray}
The first three terms in \reff{eq.proof.thm2.b}
give $a_0(y) = 1$ and $a_1(y) = 1 - y/(1-y)$,
so that $1 - a_1(y)/a_0(y) = y/(1-y) \myge 0$.
On the other hand, the final sum in \reff{eq.proof.thm2.b}
is manifestly a power series with nonnegative coefficients in $-x$ and $y$,
which proves that $(-1)^m a_m(y) \myge 0$ for all $m \ge 2$.
\qed

{\bf Remarks.}
1.  We can obtain an explicit formula for the coefficients $a_m(y)$
by inserting into \reff{eq.proof.thm2.b}
the expansion \cite[Theorem~3.3]{Andrews_76}
\be
   {1 \over (-xy;y)_{n-1}}
   \;=\;
   \sum_{k=0}^\infty \qbinom{n+k-2}{k}{y} \, (-xy)^k
   \qquad\hbox{for } n \ge 2
\ee
where the $q$-binomial coefficients are defined by
\be
   \qbinom{n}{k}{q}  \;=\;
   {(q;q)_n \over (q;q)_k \, (q;q)_{n-k}}
   \qquad\hbox{for }  0 \le k \le n  \;.
\ee
This yields
\be
   f(x,y)
   \;=\;
%%    1 + x - {xy \over 1-y} +
%%       \sum_{k=0}^\infty (-xy)^k \sum_{n=2}^\infty
%%           {(y;y)_{n+k-2} \over (y;y)_n \, (y;y)_{n-2} \, (y;y)_k} \,
%%                 (-x)^n \, y^{n^2}
   1 + x - {xy \over 1-y} +
      \sum_{n=2}^\infty \sum_{k=0}^\infty
          \qbinom{n+k-2}{k}{y} \: (-xy)^k \: {(-x)^n \, y^{n^2} \over (y;y)_n}
   \;.
\ee
Extracting the coefficient of $x^m$ for $m=n+k \ge 2$, we have
\be
   (-1)^m a_m(y)
   %% & = &
   %% \sum_{n=2}^m {(y;y)_{m-2} \over (y;y)_n \, (y;y)_{n-2} \, (y;y)_{m-n}} \,
   %%     y^{n^2 + m-n}   \\[2mm]
   \;=\;
   y^m \sum_{n=2}^m \qbinom{m-2}{m-n}{y} \, {y^{n(n-1)} \over (y;y)_n}
      \;.
 \label{eq.amy.b}
\ee
Since the $q$-binomial coefficients
are polynomials in $q$ with nonnegative integer coefficients
\cite[Theorem~3.2 or 3.6]{Andrews_76},
we see once again that $(-1)^m a_m(y) \myge 0$ for all $m \ge 2$.
We also see from \reff{eq.amy.b} that $a_m(y)$
is a rational function of the form $a_m(y) = P_m(y)/(y;y)_m$
where $P_m(y)$ is a polynomial with integer coefficients.
%% \hspace{-1mm}
%% \qed

2.  It would be interesting to seek a combinatorial interpretation
of the coefficients of $1-1/\xi_0(y)$,
analogously to what Prellberg \cite{Prellberg_11} has done for $\xi_0(y)$
[see footnotes~\ref{note_Prellberg_1}--\ref{note_Prellberg_3} above].

% \begin{lemma}
%   \label{lemma.series}
% {\bf No, this is false!!!!  It seems to be true for certain ranges
%    of the pair $(k,\ell)$, which I am trying to determine empirically.
%    But anyway it does seem to be true for $\ell=2$ and $k \ge 1$.}
% For integers $k \ge 0$ and $\ell \in \Z$,
% \be
%    S_{k,\ell}(y)
%    \;=\;
%    \sum_{n=0}^\infty  \qbinom{n+k}{k}{y} \, (y^{n+\ell+1};y)_\infty \; y^n
% \ee
% is a power series in $y$ with nonnegative coefficients.
% {\bf I believe it is also a rational function of $y$ --- with what
%    denominator???}
% \end{lemma}
% 
% 
% \proof
% ??????????????????
% \qed

\section{Proof of Theorem~\ref{thm3}}   \label{sec.proof3}

Next we prove Theorem~\ref{thm3}.
It is convenient once again to state and prove first an abstract result
\cite{Sokal_Fxy_powerseries},
and then verify the hypotheses of this abstract result in our specific case.

\begin{proposition}
  \label{prop.scrs-2.abstract}
Consider a formal power series $f(x,y)$ satisfying all the hypotheses
of Proposition~\ref{prop.scrs-1.abstract}.
Then
\be
   \xi_0(y)^{-2}
   \;\myle\;
   \biggl( {a_1(y) \over a_0(y)} \biggr)^{\! 2}  \,-\,
      2 \sum_{n=2}^\infty (-1)^n \, {a_n(y) \over a_0(y)} \,
                          \biggl( {a_0(y) \over a_1(y)} \biggr)^{\! n-2}
   \;.
 \label{eq.prop.scrs-2.abstract}
\ee
\end{proposition}

\proof
Divide both sides of \reff{eq.xi0.recursion.b} by $\xi_0(y)$
and then insert \reff{eq.xi0.recursion.b} in the first term
on the right-hand side:  we obtain
\begin{subeqnarray}
   \xi_0(y)^{-2}
   & = &
   {a_1(y) \over a_0(y)} \, \xi_0(y)^{-1} \,-\,
      \sum_{n=2}^\infty (-1)^n \, {a_n(y) \over a_0(y)} \, \xi_0(y)^{n-2}
            \\[2mm]
   & = &
   \biggl( {a_1(y) \over a_0(y)} \biggr)^{\! 2}   \,-\,
      \sum_{n=2}^\infty (-1)^n \, {a_n(y) \over a_0(y)} \,
      \biggl[ 1 + {a_1(y) \over a_0(y)} \, \xi_0(y) \biggr] \, \xi_0(y)^{n-2}
   \;.
     \nonumber \\[-4mm]
 \slabel{eq.xi0.recursion.Sminus2.b}
 \label{eq.xi0.recursion.Sminus2}
\end{subeqnarray}
Now, by hypothesis we have $(-1)^n a_n(y)/ a_0(y) \myge 0$ for all $n\ \ge 2$.
By Proposition~\ref{prop.scrs-1.abstract} we have
$\xi_0(y)^{-1} \myle a_1(y)/a_0(y) \myle 1$,
hence $\xi_0(y)^{n-2} \myge [a_0(y)/a_1(y)]^{n-2} \myge 1$ for all $n \ge 2$.
Finally, multiplying \reff{eq.xi0.recursion.b} by $\xi_0(y)$
and rearranging gives
\be
   {a_1(y) \over a_0(y)} \, \xi_0(y)
   \;=\;
   1 \,+\, \sum_{n=2}^\infty (-1)^n \, {a_n(y) \over a_0(y)} \, \xi_0(y)^n
   \;\myge\;
   1
   \;.
\ee
Inserting these facts into \reff{eq.xi0.recursion.Sminus2.b}
proves \reff{eq.prop.scrs-2.abstract}.
\qed

\proofof{Theorem~\ref{thm3}}
We again use \reff{eq.lemma.identity2}
and thus apply Proposition~\ref{prop.scrs-2.abstract} to the power series
\reff{eq.proof.thm2.b}.
While proving Theorem~\ref{thm2}
we showed that $a_0(y) = 1$, $a_1(y) = 1 - y/(1-y) \myle 1$
and $(-1)^n a_n(y) \myge 0$ for all $n \ge 2$,
so all the hypotheses of Proposition~\ref{prop.scrs-2.abstract} are satisfied.
Furthermore, from either \reff{eq.proof.thm2.b} or \reff{eq.amy.b}
it is easy to see that
\be
   (-1)^n a_n(y)   \;\myge\;   {y^{n+2} \over (1-y)(1-y^2)}
                   \;\myge\;   {y^{n+2} \over 1-y}
\ee
for all $n \ge 2$.\footnote{
   It suffices to take the term $n=2$ in
   \reff{eq.proof.thm2.b} or \reff{eq.amy.b},
   using the fact that all other terms are $\myge 0$.
}
{}From \reff{eq.prop.scrs-2.abstract} we then have
\begin{subeqnarray}
   \xi_0(y)^{-2}
   & \myle &
   \biggl( {a_1(y) \over a_0(y)} \biggr)^{\! 2}  \,-\,
      2 \sum_{n=2}^\infty (-1)^n \, {a_n(y) \over a_0(y)} \,
                              \biggl( {a_0(y) \over a_1(y)} \biggr)^{\! n-2}
        \\[2mm]
   & \myle &
   \biggl( {a_1(y) \over a_0(y)} \biggr)^{\! 2}  \,-\,
      2 \sum_{n=2}^\infty (-1)^n \, {a_n(y) \over a_0(y)}
        \\[2mm]
   & \myle &
   \biggl( {1-2y \over 1-y} \biggr)^2  \,-\,  {2 y^4 \over (1-y)^2}
        \\[2mm]
   & = &
   1 \,-\, 2y \,-\, y^2 \,-\, \sum_{n=4}^\infty (n-3) \, y^n
   \;,
\end{subeqnarray}
which proves Theorem~\ref{thm3}.
\qed

{\bf Remarks.}
1.  If we use \reff{eq.amy.b} and expand the right-hand side of
\reff{eq.prop.scrs-2.abstract}, we obtain
\be
   \xi_0(y)^{-2}  \;\myle\;
   1-2 y-y^2 \hphantom{-y^3} - y^4-2 y^5-7 y^6-18 y^7-49 y^8-130 y^9-343 y^{10}
   - \ldots
   \,,
\ee
which differs from the exact $\xi_0(y)^{-2}$ starting at order $y^8$.
The difference at order $y^8$ arises from a contribution to $\xi_0(y)$
that is proportional to $a_2(y)^2$.
The full structure of the contributions to $\xi_0(y)$ and its powers
can be read off the explicit implicit function formula \cite{Sokal_implicit}:
see \cite{Sokal_Fxy_powerseries} for details.
%% \qed

2.  It would be interesting to seek a combinatorial interpretation
of the coefficients of $1-1/\xi_0(y)^2$,
analogously to what Prellberg \cite{Prellberg_11} has done for $\xi_0(y)$
[see footnotes~\ref{note_Prellberg_1}--\ref{note_Prellberg_3} above].

\section{Identities for $\bm{R(x,y,q)}$}    \label{sec.Rxyq}

In this section we obtain some simple identities for
the three-variable Rogers--Ramanujan function \cite{Sokal_Fxy_powerseries}
\be
   R(x,y,q)
   \;=\;
   \sum\limits_{n=0}^\infty  {x^n \, y^{n(n-1)/2}  \over  (q;q)_n}
   \;.
 \label{def.Rxyq}
\ee
%% which generalize \reff{eq.lemma.identity1}/\reff{eq.lemma.identity2}
%% and reduce to them when $q=0$.
The basic principle is in fact more general,
and applies to an arbitrary power series of the form
\be
   F(x,q)  \;=\;  \sum\limits_{n=0}^\infty {a_n x^n \over (\alpha;q)_n}
   \;.
\ee

\begin{lemma}
   \label{lemma.Rxyg_general}
For arbitrary coefficients $(a_n)_{n=0}^\infty$
and an arbitrary constant $\alpha$, we have
\be
   \sum\limits_{n=0}^\infty {a_n x^n \over (\alpha;q)_n}
   \;=\;
   {1  \over  (\alpha;q)_\infty}
   \sum\limits_{\ell=0}^\infty
       {(-\alpha)^\ell \, q^{\ell(\ell-1)/2}  \over (q;q)_\ell}
   \sum\limits_{n=0}^\infty a_n \, (q^\ell x)^n
\ee
as formal power series.
\end{lemma}

\proof
Write
\be
   \sum\limits_{n=0}^\infty {a_n x^n \over (\alpha;q)_n}
   \;=\;
   \sum\limits_{n=0}^\infty a_n x^n
         \: {(\alpha q^n;q)_\infty  \over  (\alpha;q)_\infty}
\ee
and substitute Euler's second identity \reff{eq.euler2}
for $(\alpha q^n;q)_\infty$, yielding
\begin{subeqnarray}
   \sum\limits_{n=0}^\infty {a_n x^n \over (\alpha;q)_n}
   & = &
   {1  \over  (\alpha;q)_\infty}
   \sum\limits_{n=0}^\infty  a_n x^n
   \sum\limits_{\ell=0}^\infty  {(-\alpha q^n)^\ell \, q^{\ell(\ell-1)/2}
                                 \over (q;q)_\ell}
             \qquad \\[1mm]
   & = &
   {1  \over  (\alpha;q)_\infty}
   \sum\limits_{\ell=0}^\infty
       {(-\alpha)^\ell \, q^{\ell(\ell-1)/2}  \over (q;q)_\ell}
   \sum\limits_{n=0}^\infty  a_n \, (q^\ell x)^n
   \;.
\end{subeqnarray}
\qed

Specializing to $a_n = y^{n(n-1)/2}$ and $\alpha=q$,
we obtain a simple identity that expresses $R(x,y,q)$
in terms of the partial theta function:

\begin{corollary}
   \label{cor.Rxyq}
%% The function \reff{def.Rxyq} satisfies
The three-variable Rogers--Ramanujan function \reff{def.Rxyq} satisfies
\be
   %% \sum\limits_{n=0}^\infty  {x^n \, y^{n(n-1)/2}  \over  (q;q)_n}
   R(x,y,q)
   \;=\;
   {1  \over  (q;q)_\infty}
   \sum\limits_{\ell=0}^\infty
       {(-1)^\ell \, q^{\ell(\ell+1)/2}  \over (q;q)_\ell}
       \: \Theta_0(xq^\ell,y)
 \label{eq.lemma.Rxyq}
\ee
as formal power series
and as analytic functions on $(x,y,q) \in \C \times \D \times \D$.
\end{corollary}

% \proof
% Write
% \be
% %%   \sum\limits_{n=0}^\infty  {x^n \, y^{n(n-1)/2}  \over  (q;q)_n}
%    R(x,y,q)
%    \;=\;
%    \sum\limits_{n=0}^\infty x^n \, y^{n(n-1)/2}
%          \: {(q^{n+1};q)_\infty  \over  (q;q)_\infty}
% \ee
% and substitute Euler's second identity for $(q^{n+1};q)_\infty$, yielding
% \begin{subeqnarray}
%    R(x,y,q)
%    & = &
%    {1  \over  (q;q)_\infty}
%    \sum\limits_{n=0}^\infty  x^n \, y^{n(n-1)/2}
%    \sum\limits_{\ell=0}^\infty  {(-q^{n+1})^\ell \, q^{\ell(\ell-1)/2}
%                                  \over (q;q)_\ell}
%              \qquad \\[1mm]
%    & = &
%    {1  \over  (q;q)_\infty}
%    \sum\limits_{\ell=0}^\infty
%        {(-1)^\ell \, q^{\ell(\ell+1)/2}  \over (q;q)_\ell}
%    \sum\limits_{n=0}^\infty  (xq^\ell)^n \, y^{n(n-1)/2}
%              \qquad \\[1mm]
%    & = &
%    {1  \over  (q;q)_\infty}
%    \sum\limits_{\ell=0}^\infty
%        {(-1)^\ell \, q^{\ell(\ell+1)/2}  \over (q;q)_\ell}
%        \: \Theta_0(xq^\ell,y)
%    \;.
% \end{subeqnarray}
% \qed

{}From Corollary~\ref{cor.Rxyq}
we can obtain a pair of identities for $R(x,y,q)$
that generalize \reff{eq.lemma.identity1}/\reff{eq.lemma.identity2}
and reduce to them when $q=0$:

\begin{corollary}
   \label{cor.Rxyq_bis}
We have
\begin{eqnarray}
   %% \sum\limits_{n=0}^\infty  {x^n \, y^{n(n-1)/2}  \over  (q;q)_n}
   R(x,y,q)
   & = &
   {(y;y)_\infty  \over  (q;q)_\infty}
   \sum\limits_{\ell=0}^\infty
       {(-1)^\ell \, q^{\ell(\ell+1)/2}  \over (q;q)_\ell}
       \, (-x q^\ell; y)_\infty
   \sum\limits_{n=0}^\infty
       {y^n  \over (y;y)_n \, (-x q^\ell; y)_n}
       \quad \label{eq.cor.Rxyq_bis.1}  \\[3mm]
   R(x,y,q)
   & = &
   {1  \over  (q;q)_\infty}
   \sum\limits_{\ell=0}^\infty
       {(-1)^\ell \, q^{\ell(\ell+1)/2}  \over (q;q)_\ell}
       \, (-x q^\ell; y)_\infty
      \sum_{n=0}^\infty  {(-x q^\ell)^n \, y^{n^2}  \over
                          (y;y)_n \, (-x q^\ell;y)_n}
   \label{eq.cor.Rxyq_bis.2}
\end{eqnarray}
as formal power series
and as analytic functions on $(x,y,q) \in \C \times \D \times \D$.
\end{corollary}

\proof
Just substitute \reff{eq.lemma.identity1}/\reff{eq.lemma.identity2}
into \reff{eq.lemma.Rxyq}.
\qed

\smallskip

The function $\Rtilde$ defined in \reff{def.Rtilde}
is simply the rescaled version of $R$ normalized to have
$\alpha_0 = \alpha_1 = 1$:
\be
   \Rtilde(x,y,q)
   \;=\;
   R((1-q)x,y,q)
   \;=\;
   \sum\limits_{n=0}^\infty
       {x^n \, y^{n(n-1)/2}  \over
        (1+q) (1+q+q^2) \,\cdots\, (1+q+\ldots+q^{n-1})}
   \;.
\ee
Unfortunately, I do not see how to imitate the proof of
Theorem~\ref{thm1}/Proposition~\ref{prop2}
when $-1 < q < 0$ or $0 < q \le 1$.
But perhaps I am missing something.

% \bigskip
% 
% {\bf Is there a generalization \`a la Ismail--Stanton
% \be
%    R(x,y,q,b)
%    \;=\;
%    \sum\limits_{n=0}^\infty  {(bx)^n \, (-1/b;y)_n \over  (q;q)_n}
% \ee
%    ?????  (For simplicity we do not include the Chen--Liu generalization,
%    i.e.\  we consider only the case $m=1$.)
%    Does this work?????
% 
%    Indeed, there should be a generalization to
%    $\tphio(a,b;c;y,z)$ perturbed by $1/(q;q)_n$.
% 
%    And is this a bibasic hypergeometric function?????
%    And if so, are there any useful identities for it???}

\section{Discussion}  \label{sec.discussion}

The positivity results stated in Theorems~\ref{thm1}--\ref{thm3}
can be better understood by placing them
in the following general context \cite{Scott-Sokal_powerseries}:
%% In fact, much more than the positivity stated in
%% Theorems~\ref{thm1} and \ref{thm2} appears to be true.
For $\alpha \in \R \setminus \{0\}$,
let us define
% \cite{Scott-Sokal_powerseries}
the class ${\cal S}_\alpha$ to consist of those formal power series $f(y)$
with real coefficients and constant term 1 for which the series
\begin{equation}
   {f(y)^\alpha - 1  \over \alpha}
   \;=\;
   \sum_{m=1}^\infty b_{m}(\alpha) \, y^{m}
 \label{def.scrs.alpha}
\end{equation}
has all nonnegative coefficients.
The class ${\cal S}_0$ consists of those $f$
for which the formal power series
\begin{equation}
   \log f(y)
   \;=\;
   \sum_{m=1}^\infty b_{m}(0) \, y^{m}
\end{equation}
has all nonnegative coefficients.
The containment relations between the classes $\scrs_\alpha$
are given by the following fairly easy result \cite{Scott-Sokal_powerseries}:

\begin{proposition}
 % {$\!\!\!$ \bf \protect\cite{Scott-Sokal_powerseries} \ }
 %% {$\!\!\!$ \bf \protect\cite[Proposition~2.3???]{Scott-Sokal_powerseries} \ }
   \label{prop.sa.1a}
Let $\alpha,\beta \in \R$.
Then $\scrs_\alpha \subseteq \scrs_\beta$ if and only if either
\begin{itemize}
   \item[(a)]  $\alpha \le 0$ and $\beta \ge \alpha$, or
   \item[(b)]  $\alpha > 0$ and $\beta \in \{\alpha,2\alpha,3\alpha,\ldots\}$.
\end{itemize}
\nopagebreak
Moreover, the containment is strict whenever $\alpha \neq \beta$.
\end{proposition}

For the partial theta function \reff{def.G},
Theorem~\ref{thm1} states that $\xi_0 \in \scrs_1$;
Theorem~\ref{thm2} states the stronger result that $\xi_0 \in \scrs_{-1}$
(and hence that $\xi_0 \in \scrs_\alpha$ for all $\alpha \ge -1$);
and Theorem~\ref{thm3} states the yet stronger result that
$\xi_0 \in \scrs_{-2}$
(and hence that $\xi_0 \in \scrs_\alpha$ for all $\alpha \ge -2$).
%
% I~have verified Conjecture~\ref{conj_Sminus2} through order $y^{6999}$.
% -rw-r--r--  1 as2a physics  631545 Feb  7  2010 Fxy_general_x0_v2a_Gxy.out
% -rw-r--r--  1 as2a physics 7167907 Feb  6  2010 Fxy_general_x0_v2a_Gxy_out.m%
% -rw-r--r-- 1 as2a physics     1023204 Mar 18 18:02 Fxy_general_x0_v2ab2_Gxy.out
% -rw-r--r-- 1 as2a physics 39703750521 Mar  7 02:39 Fxy_general_x0_v2ab_Gxy_out_7000.m
%  and an easy use of xihat0ser[3500] from the latter file to test its -2 power
% I wonder whether a variant of the argument given in this paper
% might allow a proof of this conjecture.
%% The proof uses identity \reff{eq.lemma.identity2}
%% together with the explicit implicit function formula \cite{Sokal_implicit}.
%
This is best possible, since from
\be
   {\xi_0(y)^\alpha - 1 \over \alpha}
   \;=\;
   y \,+\, {\alpha+3 \over 2} y^2 \,+\, {(\alpha+2)(\alpha+7) \over 6} y^3
     \,+\, O(y^4)
\ee
we see immediately that $\xi_0 \notin \scrs_\alpha$ for $\alpha < -2$.

For the deformed exponential function \reff{def.F},
I conjecture that $\xi_0 \in \scrs_{-1}$
(see also \cite[Example~4.3]{Sokal_implicit}),
%% (and hence that $\xi_0 \in \scrs_\alpha$ for all $\alpha \ge -1$),
and I~have verified this through order $y^{899}$.
% -rw-r--r--  1 as2a physics  822226 Feb  4  2010 Fxy_general_x0_v2_Fxy.out
% -rw-r--r--  1 as2a physics 1788399 Feb  2  2010 Fxy_general_x0_v2_Fxy_out.m
It follows from the asymptotics of $\xi_0(y)$ as $y \uparrow 1$
\cite{Sokal_Fxy_asymptotic}
that $\xi_0 \notin \scrs_\alpha$ for $\alpha < -1$.

\begin{figure}[t]
\centering
\includegraphics[width=0.9\textwidth]{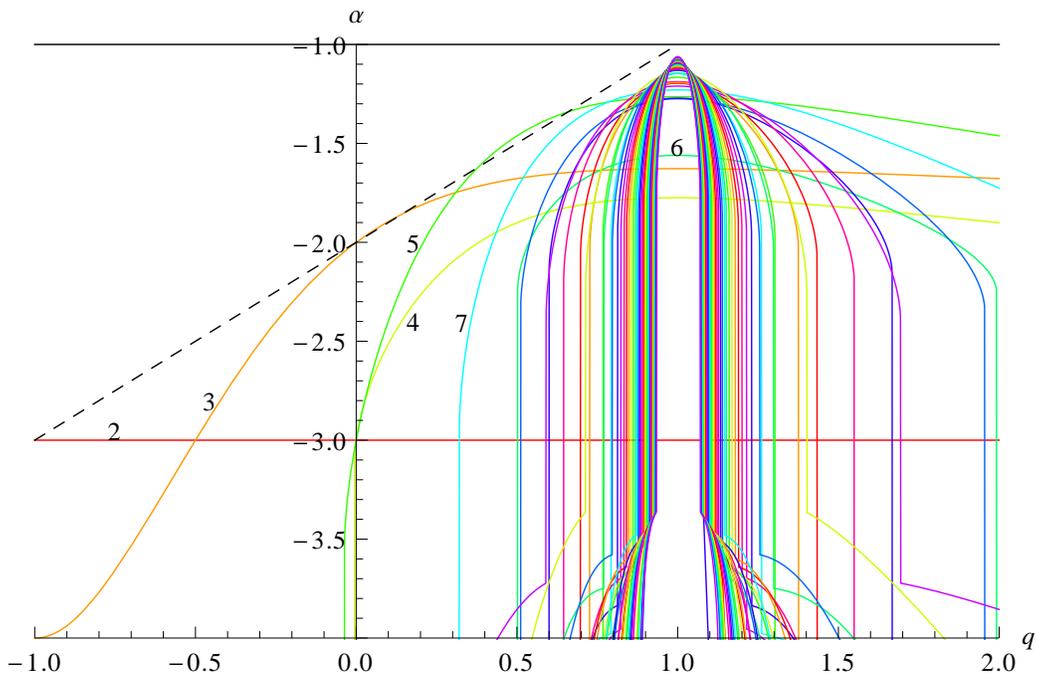}
\caption{
   Largest real root of $b_m(\alpha)$ as a function of $q$
   for $2 \le m \le 50$.
   The curves corresponding to $m \le 7$ are labeled.
   The dashed black line is $\alpha = -2+q$.
   \vspace*{5mm}
}
\label{fig.Rtilde}
\end{figure}

For the function $\Rtilde$ defined in \reff{def.Rtilde},
I conjecture that $\xi_0 \in \scrs_{-1}$ for all $q > -1$,
and I~have verified this through order $y^{349}$.
%% -rw-r--r-- 1 as2a physics 6371 Apr 13 10:48 Fxy_general_x0_v2c_Rxy_testoneoverxihat0ser.out
%% -rw-r--r-- 1 as2a physics 1472 Apr 13 06:20 Fxy_general_x0_v2c_Rxy_testoneoverxihat0ser.m
More strongly, I conjecture that for $q > -1$
there is a function $\alpha_\star(q)$
such that $\xi_0(y;q) \in \scrs_\alpha$
if and only if $\alpha \ge \alpha_\star(q)$,
and having the following properties:
\begin{itemize}
   \item[(a)]  $\alpha_\star(q) = -3$ for $-1 < q \le -1/2$.
   \item[(b)]  $\alpha_\star(q)$ is strictly increasing on $-1/2 \le q \le 1$
       and strictly decreasing on $q \ge 1$.
   \item[(c)]  $\alpha_\star(0) = -2$.
   \item[(d)]  $\alpha_\star(1) = -1$.
   \item[(e)]  $\alpha_\star(q) = \alpha_\star(1/q)$ for $q>0$.
\end{itemize}
Since
\begin{eqnarray}
%%   & &
%%   \!\!\!\!\!
   {\xi_0(y,q)^\alpha - 1 \over \alpha}
%%   \;\,=\;\,
   & \,=\, &
   {y \over 1+q} \;+\; {\alpha+3 \over 2} \: {y^2 \over (1+q)^2}  \;+\;
          O(y^3)
%%         \nonumber \\[1mm]
%%   & & \quad
%%    {(\alpha+2)(\alpha+7)(1+q^2) + (\alpha+1)(\alpha+8)q \over 6} \:
%%          {y^3 \over (1+q)^3 (1+q+q^2)}
%%     \;+\; O(y^4)
   \;,
%%         \nonumber \\[-2mm]
\end{eqnarray}
we see immediately that $\xi_0 \notin \scrs_\alpha$ for $\alpha < -3$.
Figure~\ref{fig.Rtilde} shows numerical computations
of the largest real root of $b_m(\alpha)$ [cf.\ \reff{def.scrs.alpha}],
as a function of $q \in (-1,2]$, for $2 \le m \le 50$.
%
% -rw-r--r-- 2 as2a physics 918987 May  3 04:15 plotRxy50.eps
% Created by
% -rw-r--r-- 1 as2a physics   3024 Apr 16 10:46 Fxy_general_x0_v2c_Rxy_plots.m
% -rw-r--r-- 1 as2a physics    849 May  3 04:15 Fxy_general_x0_v2c_Rxy_plots.out
% Need to run this program using "ssh -X" (otherwise it fails)
%
%% OLDER VERSION:
%% -rw-r--r--  1 as2a physics 718679 Dec 19 16:07 plotRxy40.eps
%% Created by
%%    ~/hierarchical_potts/Fxy_general_x0_v2c_Rxy_testoneoverxihat0ser.nb
%%
The upper envelope of these curves should be $\alpha_\star(q)$.
The simple conjecture $\alpha_\star(q) \le -2+q$ (shown as a dashed black line)
barely fails in the range $0 < q \ltapprox 0.145103$
because of the coefficient of $y^3$,
and in the range $0.378619 \ltapprox q \ltapprox 0.660551$
because of the coefficient of $y^5$;
but it appears to hold for $-1 < q \le 0$.
Indeed, for $-1 < q \le 0$ it appears that
$b_m(\alpha) \ge 0$ whenever $\alpha \ge -3$ and $m \neq 3$.

%
% FIGURE 1
%

% Above all, I wonder whether there might exist a generalization
% of the identity \reff{eq.identity}
% to handle the three-variable function $\Rtilde(x,y,q)$.
% If so, this might permit a proof of the conjectured positivity properties
% of $\xi_0(y,q)$.

Finally, though in this paper I have treated $\xi_0(y)$
as a formal power series,
it is not difficult to show \cite{Sokal_Fxy_asymptotic,Sokal_Fxy_hadamard},
using Rouch\'e's theorem, that $\xi_0(y)$ is in fact convergent
for $|y| < \delta_1 \approx 0.2247945929$,
where $\delta_1$ is the positive root of
$\sum_{\ell=-1}^\infty \delta^{\ell^2/2} = 2$.
(This proof applies to both $\Theta_0$ and $F$,
 and more generally to $\Rtilde$ for all $q \ge 0$.)
Then the coefficientwise positivity established in Theorem~\ref{thm1}
implies, by Pringsheim's theorem,
that the first singularity of $\xi_0(y)$ for the partial theta function
lies on the positive real axis,
namely at the point $y = y^\star_{01}$
where the leading root $x_0(y)$ collides with the next root $x_1(y)$:
this is the solution of the system
\be
   \Theta_0(x,y) \,=\, 0
   \hbox{ and }
   {\partial \Theta_0(x,y) \over \partial x} \,=\, 0
\ee
and lies at
$(x,y) = (x^\star_{01}, y^\star_{01})
 \approx (-2.3203769443,0.3092493386)$.\footnote{
   For more information concerning the real roots of the
   partial theta function and related polynomials,
   see \cite[p.~100]{Hardy_04}
   %% {\bf Also Petrovitch?????}
   \cite[pp.~330--331]{Hutchinson_23}
   \cite[vol.~1, Part~II, Problem~200, pp.~143 and 345--346,
         and vol.~2, Part~IV, Problem~176, pp.~66 and 245--246]{Polya_72}
   \cite{Kurtz_92}
   \cite[Sections~2 and 3]{Andrews_05b}
   \cite[Example~4.10]{Craven_05}
   \cite[Theorem~4]{Katkova_03}.
}
Similarly, for the deformed exponential function \reff{def.F}
it is known \cite{Morris_72,Liu_98,Langley_00}
that $\xi_0(y)$ is analytic in a complex neighborhood of the
real interval $0 < y < 1$;
therefore, if the coefficients are indeed nonnegative,
Pringsheim's theorem implies the striking fact that $\xi_0(y)$
is analytic in the whole unit disc $|y| < 1$.

\section*{Acknowledgments}

I wish to thank George Andrews, Alex Eremenko, Ira Gessel,
Christian Krattenthaler, Thomas Prellberg, Tanguy Rivoal,
Andrea Sportiello and Richard Stanley
for helpful conversations and/or correspondence.

I also wish to thank the Institut Henri Poincar\'e -- Centre Emile Borel
for hospitality during the programme on
Statistical Physics, Combinatorics and Probability
(September--December 2009), where this work was begun;
and the Laboratoire de Physique Th\'eorique
at the \'Ecole Normale Sup\'erieure (April--June 2011),
where this work was completed.

This research was supported in part
by U.S.\ National Science Foundation grant PHY--0424082.

\end{document}